\documentclass[12pt]{amsart}

\date{}
\newtheorem{lemma}{Lemma}
\newtheorem{proposition}{Proposition}

\theoremstyle{definition}

\newcommand{\IN}{{\mathbb{N}}}

\newcommand{\IZ}{{\mathbb{Z}}}

\newcommand{\IC}{{\mathbb{C}}}

\newcommand{\ind}{\mathrm{ind}}

\begin{document}
\title{Constructing subsets of a given packing index in Abelian groups}
\author{N. Lyaskovska}

\begin{abstract} By definition, the sharp packing index $\ind_P^\sharp(A)$ of a subset $A$ of an abelian group $G$ is the smallest cardinal $\kappa$ such that for any subset $B\subset G$ of size $|B|\ge\kappa$ the family $\{b+A:b\in B\}$ is not disjoint. We prove that an infinite Abelian group $G$ contains a subset $A$ with given index $\ind_P^\sharp(A)=\kappa$ if and only if one of the following conditions holds: (1) $2\le \kappa\le|G|^+$ and $k\notin \{3,4\}$; (2) 
$\kappa=3$ and $G$  is not isomorphic to $\oplus_{i\in I}
\mathbb{Z}_3$; (3) $\kappa=4$ and $G$  is not isomorphic to $\oplus_{i\in I}
\mathbb{Z}_2$ or to
 $\mathbb{Z}_4\oplus(\oplus_{i\in I} \mathbb{Z}_2)$.
\end{abstract}

\maketitle \baselineskip15pt The famous problem of optimal sphere packing traces its history back to B.Pascal and belongs to the most difficult problems of combinatorial geometry \cite{CS}. In this paper we consider an analogous problem in 
the algebraic setting. Namely, given a subset $A$ of an Abelian group $G$ we study the cardinal number
$$\ind_P(A)=\sup\big\{|B|:\mbox{$B\subset G$ and $(B-B)\cap (A-A)=\{0\}$}\big\}$$
called the {\em packing index} of $A$ in $G$. Note that the equality $(B-B)\cap (A-A)=\emptyset$ holds if and only if $(b+A)\cap (b'+A)=\emptyset$ for any distinct points $b,b'\in B$. Therefore, $\ind_P(A)$ can be thought as the maximal number of pairwise disjoint shift copies of $A$ that can be placed in the group $G$. In this situation it is natural to ask if such a maximal number always exists. In fact, this was a question of D.Dikranjan and I.Protasov who asked in \cite{DP} if for each subset $A\subset\IZ$ with $\ind_P(A)\ge\aleph_0$ there exists an infinite family of pairwise disjoint shifts of $A$. The answer to this problem turned out to be negative, see \cite{BL1}, \cite{BL2}. 
So the supremum  in the definition of $\ind_P(A)$ cannot be replaced by the maximum.

To catch the difference between $\sup$ and $\max$, let us adjust the definition of the packing index $\ind_P(X)$ and define the cardinal number
$$\ind_P^\sharp(A)=\min\{\kappa:\forall B\subset G\;\;|B|\ge\kappa\Rightarrow (B-B\cap A-A\ne\{0\})\}$$
called the {\em sharp packing index} of $A$ in $G$. In terms of the sharp packing index the question of D.~Dikranjan and I.~Protasov can be reformulated as finding a subset $A\subset \IZ$ with $\ind^\sharp_P(A)=\aleph_0$. According to \cite{BL2} (and \cite{BL1}) such a set $A$ can be found in each infinite (abelian) group $G$. Having in mind this result, I.Protasov asked in a private conversation if for any non-zero cardinal $\kappa\le|G|$ there is a set $A\subset G$ with $\ind_P(A)=\kappa$. In this paper we answer this question affirmatively (with three exceptions). Firstly, we treat a similar question for the sharp packing index because its value  completely determines the value of $\ind_P(A)$:
$$\ind_P(A)=\sup\{\kappa:\kappa<\ind_P^\sharp(A)\}.$$

Our principal result is
\smallskip

\noindent{\bf Main Theorem.} {\em An  infinite Abelian
group $G$ contains a subset $A\subset G$ with sharp packing index
$ind_P^\sharp(A)=\kappa$ if and only if one of the following conditions
holds:
\begin{enumerate}
\item [1)] $2\leq \kappa \leq |G|^+$ and $\kappa\not \in \{3,4\}$.
\item [2)] $\kappa=3$ and $G$  is not isomorphic to $\oplus_{i\in I}
\mathbb{Z}_3$.
\item [3)] $\kappa=4$ and $G$  is not isomorphic to $\oplus_{i\in I}
\mathbb{Z}_2$ or to
 $\mathbb{Z}_4\oplus(\oplus_{i\in I} \mathbb{Z}_2)$.
\end{enumerate}
}
\smallskip

Using the relation between the packing  and sharp packing indices, we can derive from the above theorem an analogous characterization of possible values of the packing index.

\smallskip

\noindent{\bf Corollary.} {\em An  infinite Abelian
group $G$ contains a subset $A\subset G$ with packing index
$ind_P(A)=\kappa$ if and only if one of the following conditions
holds:
\begin{enumerate}
\item [1)] $1\leq \kappa \leq |G|$ and $\kappa\not \in \{2,3\}$.
\item [2)] $\kappa=2$ and $G$  is not isomorphic to $\oplus_{i\in I}
\mathbb{Z}_3$.
\item [3)] $\kappa=3$ and $G$  is not isomorphic to $\oplus_{i\in I}
\mathbb{Z}_2$ or to
 $\mathbb{Z}_4\oplus(\oplus_{i\in I} \mathbb{Z}_2)$.
\end{enumerate}}
\smallskip

\section{Preliminaries}

In the proof of Main Theorem we shall exploit a combinatorial lemma proved in this section. For a set $A$ by $[A]^2=\{B\subset A:|B|=2\}$ we denote the family of all two-element subsets of $A$.

We shall say that a map $f: [A]^2\mapsto [B]^2$ 

\begin{itemize}
\item is {\em separately
injective} if for any $a\in A$ the map $f_a: x\mapsto f(\{x,a\})$
is injective;
\item {\em preserves intersections} if for any $a_0,a_1,a_2\in A$ the
intersection $f(\{a_0,a_1\})\cap f(\{a_0,a_2\})$ is not empty.
\end{itemize}

\begin{lemma}\label{size} If  $|A|\ge 5$  and a map
$f:[A]^2\mapsto [B]^2$ is separately injective and preserves
intersections, then $|A|\leq |B|$.
\end{lemma}

\begin{proof}
Fix any point $a_0\in A$ and consider the family
$\big\{f(\{a,a_0\}):a\in A\backslash \{a_0\}\big\}$.
 Since $f$ preserves intersections we have that
  $f(\{a,a_0\})\cap f(\{a^{'},a_0\})\ne \emptyset$ for any
  distinct $ a,a^{'}\in A$.
Using the separately injective of $f$ and the inequality $|A|\ge5$ we
can prove that the intersection 
 $\bigcap_{a\in A\backslash \{a_0\}}f(\{a,a_0\})$ is not empty and hence contains some element $b_0$. Thus we obtain that $f:\{a,a_0\}\mapsto \{b,b_0\}.$ And since $f$ is
separately injective we obtain an injective map from  $A\backslash
\{a_0\}$ into $B\backslash \{b_0\}$ implying the desired inequality $|A|\leq|B|$.
\end{proof}

We shall also need one structure property of Abelian groups. By $\IZ$ we denote the additive group of integer numbers and by
$$\IZ(p^\infty)=\{z\in \IC:\exists n\in\IN\mbox{ with }z^{p^k}=1\}$$ the quasicyclic $p$-group for a prime number $p$.

\begin{proposition}\label{structure} Each infinite Abelian group $G$ contains an infinite subgroup isomorphic to $\IZ$, $\IZ(p^\infty)$ or the direct sum of finite cyclic groups.
\end{proposition}

\begin{proof} If  $G$ contains an element $g$ of infinite
order, then it generates a cyclic subgroup isomorphic to $\IZ$. 
Otherwise, $H$ is a torsion group and by Theorem 8.4 \cite{Fu} decomposes into the direct sum $G=\oplus_p A_p$ of $p$-groups
$A_p$. If each group $A_p$ is finite, then $G$ contains an infinite direct product of finite cyclic group. If for some prime number $p$ the $p$-group $A_p$ is infinite, then there are two cases. Either $A_p$ contains a copy of the quasicyclic $p$-group $\IZ(p^\infty)$ or else each element of $A_p$ has finite height. In the latter case, take any infinite countable subgroup $H\subset A_p$ and apply Theorem 17.3 of \cite{Fu} to conclude that $H$ is the direct sum of finite cyclic groups.
\end{proof}

\section{The proof of the ``only if'' part of  Main Theorem}

The proof of the ``only if'' part of Main Theorem is divided into two lemmas.

\begin{lemma} If a group $G$ contains a subset $A\subset G$ with $
ind_p^\sharp(A)= 3$ (which is equivalent to $\ind_P(A)=2$), then $G$ is not  isomorphic to the direct sum
$\oplus_{i\in I} \mathbb{Z}_3$.

\end{lemma}

\begin{proof}

On the contrary suppose that $G$ is   isomorphic to the direct sum
$\oplus_{i\in I} \mathbb{Z}_3$ and take a subset $A$ of $G$
with $\ind_P^\sharp(A)=3$. The latter is equivalent to $\ind_P(A)=2$
which means that there is a subset $B_2\subset G$ of size $2$ such
that the family $\{b+A: b\in B_2\}$ is disjoint. Note that for
every $b^{'}\in G$ the family $\{b+A: b\in b^{'}+B_2\}$ is 
disjoint too. So without loss of generality we can assume that
$B_2=\{0,b_1\}$. The family $\{b+A: b\in B_2\}$ is disjoint and hence
$$A\cap (b_1+A)=\emptyset.$$

Adding to both sides $b_1$ and $2b_1$ we get
$$(b_1+A)\cap (2b_1+A)=\emptyset;$$
$$(2b_1+A)\cap (3b_1+A)=\emptyset.$$

 Since $G$ is isomorphic to the direct sum $\oplus_{i\in I}
\mathbb{Z}_3$ we get $3b_1=0$. Thus we conclude that $\big \{ b+A: b\in
\{0,b_1,2b_1\}\big\}$ is disjoint and so $\ind_P(a)>2$ and
$ind_p^\sharp(A)>3$, which contradicts our assumption.
\end{proof}

\begin{lemma} If a group $G$ contains a subset $A\subset G$ with  $
ind_p^\sharp(A)= 4$ (which is equivalent to $\ind_P(A)=3$), then $G$ can not be isomorphic neither  to the
direct sum $\oplus_{i\in I} \mathbb{Z}_2$ nor to
$\mathbb{Z}_4\oplus(\oplus_{i\in I} \mathbb{Z}_2)$.
\end{lemma}

\begin{proof}

Conversely suppose that $G$ is  isomorphic to $\oplus_{i\in I}
\mathbb{Z}_2$ or to $\mathbb{Z}_4\oplus(\oplus_{i\in I}
\mathbb{Z}_2)$ and there exists a subset $A$ of $G$ with
$\ind_P^\sharp(A)=4$. This is  equivalent to $\ind_P(A)=3$ and  from the
definition we get that there is a three-element subset $B_3\subset
G$ such that the family $\{b+A: b\in B_3\}$ is disjoint. Note that
for any $b^{'}\in G$ the family $\{b+A: b\in b^{'}+B_3\}$ is 
disjoint too. So, without loss of generality we can assume that
$B_3=\{0,b_1,b_2\}$. Since the family $\{b+A: b\in B_3\}$ is
disjoint we conclude that

\hspace{40mm}  (1) \hspace{10mm} $A\cap (b_1+A)=\emptyset;$

\hspace{40mm} (2)  \hspace{10mm} $A\cap ( b_2+A)=\emptyset;$

\hspace{40mm}  (3) \hspace{10mm} $ (b_1+A)\cap (b_2+A)=\emptyset.$

We consider three  cases.

\vspace{2mm} \textbf{Case 1.}  Suppose one of the elements
$b_1,b_2$ is of order $2$. Let it be $b_1$. Then  $2b_1=0$ and

\hspace{30mm}  $(2)+b_1:$ \hspace{10mm}$(b_1+A)\cap
(b_1+b_2+A)=\emptyset;$

\hspace{30mm}  $(3)+b_1:$ \hspace{10mm} $A\cap (b_1+
b_2+A)=\emptyset.$

\hspace{30mm}  $(1)+b_2:$ \hspace{10mm} $(b_2+ A)\cap
(b_2+b_1+A)=\emptyset.$

 Thus we get that the family $\big\{b+A: b\in
\{0,b_1,b_2,b_1+b_2\}\big\}$ is disjoint and hence $\ind_P(A)>3$
and $\ind_p^\sharp(A)>4$, which contradicts our assumption. Thus we
complete the proof of the {Case 1.}
\smallskip

Next we consider two cases where both $b_1$ and $b_2$ are of order $4$.
In this case the  group $G$ is isomorphic to $\mathbb{Z}_4\oplus
(\oplus_{i\in I} \mathbb{Z}_2).$ Therefore there are two
possibilities: $b_1=(g,x), b_2=(g,y)$ or $b_1=(g,x), b_2=(-g,y)$
where $x,y \in \oplus_{i\in I}\mathbb{Z}_2$ and  $g\in
\mathbb{Z}_4$ is of order $4$.
\smallskip

\vspace{2mm} \textbf{Case 2.} Suppose $b_1=(g,x), b_2=(g,y)$ where
$x,y \in \oplus_{i\in I}\mathbb{Z}_2$ and  $g\in \mathbb{Z}_4$ is
of order $4$.

Recall that $B_3= \big\{(0,0), (g,x),(g,y)\big\}$ and consider the set $B_4= \big\{(0,0), (g,x),(g,y),(0,x+y)\big\}$. We claim that
the family $\{b+A:b\in B_4\}$ is disjoint. Indeed, since $\{b+A:b\in B_3\}$ is disjoint we have:

\hspace{40mm}  $(1)$ \hspace{10mm} $ A\cap ((g,x)+A)=\emptyset;$

\hspace{40mm}  $(2)$ \hspace{10mm} $A\cap( (g,y)+A)=\emptyset;$

\hspace{40mm}  $(3 )$ \hspace{10mm} $ ((g,x)+A)\cap
((g,y)+A)=\emptyset.$

Then

\hspace{20mm}  $(3)+(3g,y)$\hspace{5mm} $:$ \hspace{10mm}$
((0,x+y)+A)\cap A=\emptyset;$

\hspace{20mm}  $(2)+ (0,x+y):$ \hspace{10mm}$ ((0,x+y)+A)\cap
((g,x)+A)=\emptyset;$

\hspace{20mm}  $(1)+(0,x+y):$ \hspace{10mm}$((0,x+y)+A)\cap
((g,y)+A)=\emptyset.$

Hence, the family $\{b+A:b\in B_4\}$ is disjoint which implies $\ind_P(A)\ge 3$ and  $\ind_P^\sharp(A)\ge 4$, a contradiction
with the assumption.
\smallskip

{\bf Case 3.} Suppose $b_1=(g,x), b_2=(-g,y)$ where $x,y
\in\oplus_{i\in I}\mathbb{Z}_2$ and  $g\in \mathbb{Z}_4$ is of order $4$.
\smallskip

In this case $B_3= \{(0,0), (g,x),(-g,y)\}.$ 
Put $B_4= \{(0,0), (g,x),(-g,y),(2g,x+y)\}$. We claim that the
family $\{b+A:b\in B_4\}$ is disjoint. Indeed,
since $\{b+A:b\in B_3\}$ is disjoint we have:

\hspace{40mm}  $(1)$ \hspace{10mm} $ A\cap ((g,x)+A)=\emptyset;$

\hspace{40mm}  $(2)$ \hspace{10mm} $A\cap( (-g,y)+A)=\emptyset;$

\hspace{40mm}  $(3 )$ \hspace{10mm} $ ((g,x)+A)\cap
((-g,y)+A)=\emptyset.$

Then

\hspace{20mm}  $(3)+(g,y)$\hspace{10mm} $:$ \hspace{10mm}$
((2g,x+y)+A)\cap A=\emptyset;$

\hspace{20mm}  $(2)+ (2g,x+y):$ \hspace{10mm}$ ((2g,x+y)+A)\cap
((g,x)+A)=\emptyset;$

\hspace{20mm}  $(1)+(2g,x+y):$ \hspace{10mm}$((2g,x+y)+A)\cap
((-g,y)+A)=\emptyset.$

Hence the family $\{b+A:b\in B_4\}$ is disjoint
and thus $\ind_P(A)> 3$ and $\ind_P^\sharp(A)> 4$, which contradicts our
assumption.
\end{proof}

Thus if $G$ contains a subset $A\subset G$ with
$\ind_P^+(A)=\kappa$ then one of the condition 1)-3) holds.

\section{The proof of the ``if'' part of Main Theorem}

To prove the ``if'' part of the Main Theorem, given a cardinal $\kappa$ satisfying one
of the conditions 1)--3) we shall construct a subset $A$ with
$\ind_P^\sharp(A)=\kappa$. First we shall construct a subset $A_\kappa$
assuming that we have in disposal
 an auxiliary subset $\mathbb{B}_\kappa$ with some properties. Next, a 
subset $\mathbb{B}_\kappa$ wil the desired properties will be constructed in each group.

\begin{proposition} An infinite Abelian group $G$ contains a subset
$A_\kappa$ with $\ind_P^\sharp (A_\kappa)=\kappa$ if there exists  a
subset $\mathbb{B}_\kappa=-\mathbb{B}_\kappa$ of $G$ with the
following properties:
\begin{enumerate}
\item [($1_\kappa$)]for every cardinal $\alpha < \kappa$ there is
a  subset $B_{\alpha}$ of size $|B_{\alpha}|= \alpha$ such that
$B_{\alpha}-B_{\alpha} \subset \mathbb{B}_\kappa$;
\item [($2_\kappa$)]
$B_{\kappa}-B_{\kappa} \not \subset \mathbb{B}_\kappa $ for any
subset $B_{\kappa} \subset G$ of size $\kappa$;
\item [($3_\kappa$)]$F+\mathbb{B}_\kappa\ne G$ for any subset $F\subset
G$ of size $|F|<|G|$.
\end{enumerate}

\end{proposition}

By $|A|$ we denote the cardinality of a set $A$.

\begin{proof}
Let $\mathbb{B}_\kappa^\circ=\mathbb{B}_\kappa\setminus\{0\}$. We
shall construct a subset $A_\kappa\subset G$ such that $
(\mathbb{B}_\kappa^\circ +A_\kappa) \cap A_\kappa =\emptyset$.
Moreover, the subset $A_\kappa$ will be constructed so that $
G\backslash \mathbb{B}_\kappa^\circ \subset A_\kappa-A_\kappa$.

Let $\lambda=|G\backslash \mathbb{B}_\kappa^\circ |$ and
$G\backslash \mathbb{B}_\kappa^\circ=\{ g_\alpha : \alpha <
\lambda\}$ be an enumeration of $G\backslash
\mathbb{B}_\kappa^\circ $ by ordinals $ \alpha<\lambda$.

We put  $A_\kappa=\bigcup_{\alpha<\lambda}\{a_{\alpha},g_{\alpha}+
a_{\alpha}\}$, where a sequence $(a_\alpha)_{\alpha<\lambda}$ is
to be defined later. This clearly forces that $ G\backslash
\mathbb{B}_\kappa^\circ \subset A_\kappa-A_\kappa$.

The task is now to find a sequence $(a_\alpha)_{\alpha<\lambda}$
such that $ (\mathbb{B}_\kappa^\circ +A_\kappa) \cap A_\kappa
=\emptyset$. We define this sequence by induction.

We start with $a_0=0$. Assuming that for some $\alpha$ the points $a_\beta,
\beta<\alpha$, have been constructed, put $F_{\alpha}=\{
a_\beta,g_\beta +a_\beta:\beta < \alpha \} $.

According to the property $(3_\kappa)$ of the set
$\mathbb{B}_\kappa$ we can pick
 a point $a_\alpha \in G$ so that
$$ a_\alpha \notin F_\alpha +\mathbb{B}_\kappa\cup F_\alpha
-g_{\alpha}+\mathbb{B}_\kappa .$$

This gives $ (\mathbb{B}_\kappa^\circ +A_\kappa) \cap A_\kappa
=\emptyset$.

It remains to show that $A_\kappa$ satisfies the conclusion of the
theorem.

According to the property $(1_\kappa)$ of the set
$\mathbb{B}_\kappa$ for any cardinal $\alpha<\kappa$ there is
$B_\alpha$ such that $B_\alpha-B_\alpha \subset
\mathbb{B}_\kappa$.  From the fact that $\mathbb{B}_\kappa^\circ
+A_\kappa\cap A_\kappa=\emptyset$  we conclude that $b-b^{'}
+A_\kappa\bigcap A_\kappa=\emptyset$ for all distinct $b,b^{'}\in
B_\alpha.$  Thus for any cardinal $\alpha<\kappa$ there is
$B_\alpha$ such that the family  $ \{ b+A_\kappa : b\in B_\alpha
\}$ is disjoint and so $\ind_P^\sharp(A_\kappa)\ge \kappa $.

Let us show that $\ind_P^\sharp(A_\kappa)=\kappa$.
According to the property $(2_\kappa)$, for any subset
$B_{\kappa}\subset G$ of size $ \kappa$ there are $b,b^{'}\in
B_{\kappa}$ such that $b-b^{'}\notin \mathbb{B}_\kappa$. Therefore
$b-b^{'}\in G\backslash \mathbb{B}^\circ_\kappa \subset
A_\kappa-A_\kappa.$ Hence $ b+A_\kappa\bigcap b^{'}+A_\kappa\ne
\emptyset$, which yields $\ind_P^\sharp(A_\kappa)\leq \kappa$. Combining the two inequalities, we get $\ind_P^\sharp(A_\kappa)= \kappa$.
\end{proof}

The proof of the Main Theorem will be completed as soon as we construct a subset $\mathbb{B}_\kappa$ with properties
$(1_\kappa)-(3_\kappa)$. This will be done in the following five
lemmas.

\begin{lemma} Let  $\kappa=3$  and $G$ be an  infinite Abelian
group which is not isomorphic to the direct sum $\oplus_{i\in I}
\mathbb{Z}_3$. Then $G$ contains a subset $\mathbb{B}_3$ with the
properties $(1_3 )-(3_3).$

\end{lemma}

\begin{proof}
Pick any nonzero point $g\in G$ whose order is not equal to 3 and consider the set
$\mathbb{B}_3=B_2-B_2=\{0,\pm g\}$ where $B_2=\{0,g\}$. It is clear
that $\mathbb{B}_3$ has the properties $(1_3 )$, $(3_3)$. So it is
enough to show that $\mathbb{B}_3$ satisfies the property $(2_3)$.
Note that if $2g=0$ then $\mathbb{B}_3=\{0,g\}$ is a subgroup of
$G$ and hence has the property $(2_3)$.

 So we assume that $2g\ne 0$ which yields that $\mathbb{B}_3=\{0,g,-g\}$ contains three elements. 
To  prove that $\mathbb{B}_3$ has property $(2_3)$ fix some subset
$B_3\subset G$ of size $3$ and pick any point $b_0\in B_3$. If
there is $b\in B_3$ with
$$b-b_0\not \in \mathbb{B}_3=\{0,g,-g\}$$
then there is nothing to prove. Otherwise we have that
$$B_3-b_0\subset \mathbb{B}_3.$$

Since $|B_3|=3$ there are $b, b^{'}\in B_3$ such that $b-b_0=g; \,
b^{'}-b_0=-g.$ Hence we get  $b-b^{'}=2g$. From the choice of
element $g$ we get that $2g\not \in \mathbb{B}_3$. Hence
$b_2-b_3\not \in \mathbb{B}_3$ and $\mathbb{B}$ has the property
$(2_3)$ which completes the proof of the lemma.
\end{proof}

\begin{lemma} Let  $\kappa=4$ and $G$ be an  infinite Abelian
which is not isomorphic to $\oplus_{i\in I} \mathbb{Z}_2$ or to
 $\mathbb{Z}_4\oplus (\oplus_{i\in I} \mathbb{Z}_2)$. Then $G$ contains
a subset $\mathbb{B}_4$ with properties $(1_4)-(3_4)$.
\end{lemma}

\begin{proof}
We consider three cases.
\smallskip

\vspace{2mm} \textbf{Case 1.} Suppose a group  $G$ contains an
element $g$ with order $>5$.
\smallskip

Put $\mathbb{B}_4=B_3-B_3=\{0,\pm g,\pm 2g\}$ where $B_3=
\{0,g,-g\}$. It is easily to check that $\mathbb{B}_4$ has the
properties $(1_4),(3_4).$ We claim that $\mathbb{B}_4$ satisfies the property $ (2_4).$

To derive a contradiction, suppose that there is a subset $B_4\subset G$ of size $|B_4|=4$ such that $B_4-B_4\subset \mathbb{B}_4=\{0,g,-g,2g,-2g\}$.

 Fix some element $b_0\in B_4$. Since $B_4-b_0\subset
 \mathbb{B}_4$ there are $b,b^{'}\in B_4$ such that

 $b-b_0=-g;b^{'}-b_0=2g$ or  $b-b_0=g;b^{'}-b_0=-2g.$

Then  $b^{'}-b=3g$ or $b^{'}-b=-3g$.

Note that since the order of $g$ is greater than 5, neither $3g\in \mathbb{B}_4$ no
$-3g\in \mathbb{B}_4$. Thus we get $b^{'}-b\not \in \mathbb{B}_4$,
a contradiction with the assumption.
Hence  $\mathbb{B}_4$ satisfies the property $ (2_4)$ and we
complete the proof of Case 1.
\smallskip

{\bf Case 2.} Assume that $G$ contains no element of order greater than 5. Then $G$ is the direct sum of cyclic groups according to Theorem 17.2 of \cite{Fu}. More precisely, $G$ is
isomorphic either to $(\oplus_{i\in I} \mathbb{Z}_2)\oplus (\oplus_{j\in
J}\mathbb{Z}_4 )$ or to $\oplus_{i\in I} \mathbb{Z}_3$ or to
$\oplus_{i\in I} \mathbb{Z}_5$. Since $G$ is not isomorphic to $\oplus_{i\in I}\IZ_2$ or $\IZ_4\oplus\oplus_{i\in I}\IZ_2$, we have to consider the
following two cases: $G$ contains a subgroup isomorphic to $\IZ_3$ and $G$ is contains a subgroup isomorphic to $\mathbb{Z}_i \oplus
\mathbb{Z}_j\oplus H $ for some $4\leq i,j\leq 5$.
\smallskip

\vspace{2mm} \textbf{Case 2a.} Suppose that  $G$ contains a subgroup $H$ isomorphic to $\IZ_3$.
\smallskip

In this sace we put $\mathbb{B}_4=H$ and see that $\mathbb{B}_4$ has the properties $(1_4)-(3_4)$.
\smallskip

\vspace{2mm} \textbf{Case 2b.} Suppose $G$ contains a subgroup isomorphic to the
direct sum of $\mathbb{Z}_i\oplus\mathbb{Z}_j\oplus H$ for some
$4\leq i,j\leq 5$.
\smallskip

We shall identify  $\mathbb{Z}_i\oplus\mathbb{Z}_j$ with a subgroup of $G$ and shall find a subset $\mathbb{B}_4\subset
\mathbb{Z}_i\oplus\mathbb{Z}_j$ with the properties $(1_4)-(3_4)$. Obviously
$\mathbb{B}_4$ has the same properties in the whole group $G$.

Put $\mathbb{B}_4=B_3-B_3$ where $B_3= \{(0,0),(g_1,0),(0,g_2)\}$.
It is clear that  $\mathbb{B}_4$ has the properties $(1_4),(3_4)$. We
claim that  $\mathbb{B}_4$ has property $(2_4).$
Indeed, assuming the converse, we would find a subset $B_4\subset G$ of size $|B_4|=4$ with $B_4-B_4 \subset B_3-B_3$.

Fix any point $b_0\in B_4$. Then
$$B_4-b_0\subset \mathbb{B}_4=\{(0,0),(g_1,0),(0,g_2),
(-g_1,0),(0,-g_2),(g_1,-g_2),(-g_1,g_2)\}.$$

Let us show that $(g_1,0)\not\in B_4-b_0$. Since the elements $g_1$ and $g_2$ have order $\ge 4$,
$$
\begin{aligned}
(g_1,0)-(-g_1,0)&\not \in \mathbb{B}_4;\\
(g_1,0)-(-g_1,g_2)&\not \in \mathbb{B}_4;\\
(g_1,0)-(0,-g_2)&\not \in \mathbb{B}_4.
\end{aligned}
$$

Thus if there is $b\in
B_4$ with $b-b_0=(g_1,0)$ then
$$B_4-b_0\subset \mathbb{B}_4=\{(0,0),(g_1,0),(0,g_2),
(g_1,-g_2)\}.$$

From the above and the fact that $|B_4|=4$ we get that there are
$b_1,b_2\in B_4$ such that $b_1-b_0=(0,g_2)$ and $b_2-b_0=(g_1,-g_2)$.
Hence $b_2-b_1=(g_1,-2g_2)\not \in \mathbb{B}_4$, a contradiction
with the assumption that $B_4-B_4\subset \mathbb{B}_4$. So, we 
conclude that $(g_1,0)\not\in B_4-b_0$.

In the same manner we can show that none of the elements $
(0,g_2),(-g_1,0),(0,-g_2)$ belong to $B_4-b_0$, which contradicts the fact that $B_4-B_4 \subset \mathbb{B}_4$. This completes
the proof of Lemma.
\end{proof}

\begin{lemma} If $\kappa >4$ is a finite cardinal, then each infinite Abelian group $G$contains a subset $\mathbb{B}_\kappa$ with the properties $(1_\kappa)-(3_\kappa)$.
\end{lemma}

\begin{proof}
It is easy to check that each subset $\mathbb{B}_\kappa$ with the
properties $(1_\kappa)-(3_\kappa)$ in a subgroup $H\subset G$  has
these properties in the whole group $G$. This observation combined with Proposition~\ref{structure} reduces the problem to
constructing a set $\mathbb{B}_\kappa$ in the groups $\IZ$,
$\IZ(p^{\infty})$ or the direct sum of finite cyclic groups.
This will be done separately in the following three cases.
\smallskip

\vspace{2mm} \textbf{Case 1.} We construct a subset $\mathbb B_k$ in the group $\IZ$.  
\smallskip

In this case put
$\mathbb{B}_\kappa=B_{\kappa -1}-B_{\kappa -1}$ where $B_{\kappa
-1}=\{i : 1\leq i\leq \kappa-1\}$.  It is easy to check that
$\mathbb{B}_\kappa$ has property $(1_\kappa)-(3_\kappa)$ in $\IZ$.

\smallskip

\vspace{2mm} \textbf{Case 2.} We construct a subset
$\mathbb{B}_\kappa$ in the quasicyclic $p$-group $\IZ(p^\infty)$.
\smallskip

Choose $n$ such that $z^{p^n}\in \{e^{i\phi} :\frac{2\pi
}{\kappa}<\phi<\frac{2\pi }{\kappa-1}\}$.
 Then put $\mathbb{B}_\kappa=B_{\kappa -1}-B_{\kappa -1}$ where
$$B_{\kappa -1}=\{e^{i\phi} : \phi =\frac{2\pi l}{p^n},\; 1\leq l\leq
\kappa-1\}.$$
It is easy to check that $\mathbb{B}_\kappa$ has the properties
$(1_\kappa)-(3_\kappa)$ in $\IZ(p^\infty)$.
\smallskip

\vspace{2mm} \textbf{Case 3.} We construct a subset
$\mathbb{B}_\kappa$ in the direct sum of cyclic groups
$\oplus_{i\in \omega} \langle g_i \rangle$.
\smallskip

Put $\mathbb{B}_\kappa=B_{\kappa -1}-B_{\kappa -1}$ where
$B_{\kappa -1}=\{g_i : 1\leq i\leq \kappa-1\}$. Obviously
$\mathbb{B}_\kappa$ has properties $(1_\kappa),(3_\kappa)$.
We claim that  $\mathbb{B}_\kappa$ has property $(2_\kappa).$
To obtain a contradiction assume that there exists a subset $B_\kappa
\subset G$ with size $|B_\kappa|=\kappa$ such that
$$B_\kappa-B_\kappa \subset
B_{\kappa-1}- B_{\kappa-1}.$$

Consider the sets $S=\{i: 1\leq i\leq \kappa-1\}$ and $F=\{i:
1\leq i\leq \kappa\}$.  We can enumerate the sets $B_{\kappa -1}$
and $B_\kappa$ as  $B_{\kappa-1}=\{g_i : i\in S\}$ and
$B_\kappa=\{b_i : i\in F\}$.

Since $B_\kappa-B_\kappa \subset B_{\kappa-1}- B_{\kappa-1}$ we
can define a map $f: [F]^2\mapsto [S]^2$ assigning to each pair
$\{i,j\}\in [F]^2$ a unique pair $\{k,l\}\in [S]^2$ such that
 $b_i-b_j=\pm(g_k-g_l)$.
A desired contradiction will follow from Lemma~\ref{size} as soon as we check that $f$ is separately injective and preserves
intersections.
\smallskip

\vspace{2mm}\textbf{Claim 1.} The  map $f$ preserves
intersections.
\smallskip

 To derive a contradiction, suppose that there are distinct $i,i^{'}\in F$
 and $j\in F$ such that
 $$f(\{i,j\})\cap f(\{i^{'},j\})=\emptyset.$$
Then
$$b_i-b_j=g_k-g_l\mbox{ and }b_{i^{'}}-b_j=g_n-g_m$$ where $k,l,n,m$ are pairwise distinct. 

Hence $ b_{i^{'}}-b_i=g_n-g_m-g_k+g_l \not \in
B_{\kappa-1}-B_{\kappa-1}$, which contradicts the assumption that $B_\kappa-B_\kappa
\subset B_{\kappa-1}- B_{\kappa-1}.$

\smallskip

\vspace{2mm}\textbf{Claim 2.} The map $f$ is separately injective.
\smallskip

To derive a contradiction, suppose that there are distinct $i,i^{'}\in F$
 and $j\in F$ such that
 $$f(\{i,j\})= f(\{i^{'},j\})=\{k,l\}.$$

Since $b_i,b_{i^{'}}$ are distinct we get
$$b_i-b_j=g_k-g_l\mbox{ and }b_{i^{'}}-b_j=g_l-g_k$$and thus
$b_{i^{'}}-b_i=2(g_l-g_k)\ne 0$. Note that $2(g_l-g_k)\in
B_{\kappa-1}- B_{\kappa-1}$ iff $2(g_l-g_k)=g_k-g_l.$ Thus we get
that $3g_k=0$ and $3g_l=0.$

Since $|F|=\kappa>4$ we can chose $r\in F\backslash\{i,i^{'},j\}$.
The map $f$ preserves intersections so $f(\{r,j\})\cap \{k,l\}\ne
\emptyset .$ Also note that $f(\{r,j\})\cap \{k,l\}\ne \{k,l\}$
otherwise $b_r=b_i$ or $b_r=b_{i^{'}}.$ So without loss of
generality we can assume that $f(\{r,j\})\cap \{k,l\}=\{k\} .$

Hence  $b_r-b_j=g_s-g_k$ or $b_r-b_j=g_k-g_s$ for some $s$.

Consequently, $b_r-b_i=g_s-2g_k+g_l$ or $b_r-b_{i^{'}}=2g_k-g_s-g_l.$
 
Note that $2g_k\ne 0$ since  $3g_k=0$.

So we get $b_r-b_i=g_s-2g_k+g_l\not \in B_{\kappa-1}-B_{\kappa-1}$
or $b_r-b_{i^{'}}=2g_k-g_s-g_l\not \in B_{\kappa-1}-B_{\kappa-1}.$

 This contradicts the
assumption that $B_\kappa-B_\kappa \subset B_{\kappa-1}-
B_{\kappa-1}.$
\end{proof}

\begin{lemma} Let $\kappa$ be an infinite not limit cardinal with
$\kappa \leq |G|$ where $G$ is an infinite Abelian group. Then
there exists a subset $\mathbb{B}_\kappa$ with the properties
$(1_\kappa)-(3_\kappa)$.
\end{lemma}

\begin{proof}
Since $\kappa$ is infinite not limit cardinal there exists
cardinal  $\alpha$ such that $\kappa=\alpha^+$. Put
$\mathbb{B}_\kappa=B_\alpha-B_\alpha$ where $B_\alpha$ is any
subset of $G$ with size $|B_\alpha|=\alpha$. Obviously
$\mathbb{B}_\kappa$ satisfies property $(1_\kappa).$

Since $|\mathbb{B}_\kappa|=\alpha$ and
$|B_\kappa-B_\kappa|=\kappa=\alpha^+$ for any subset
$B_\kappa\subset G$ of size $\kappa$ we get $B_\kappa-B_\kappa\not
\subset \mathbb{B}_\kappa.$ Therefore $\mathbb{B}_\kappa$ has
property $(2_\kappa)$.

The last property $(3_\kappa)$ follows from the fact $|F|+|\mathbb B_\kappa|\le |F|\cdot |\mathbb B_\kappa|<|G|$.
\end{proof}

\begin{lemma} Let $\kappa$ be a limit cardinal and $G$
be an infinite Abelian group with $\kappa\leq |G|$. Then  there
exists a subset $\mathbb{B}_\kappa\subset G$ with the properties
$(1_\kappa)-(3_\kappa)$.
\end{lemma}

\begin{proof}
 Note that it is enough to show that each group $G$ of size
 $\kappa$ contains a
subset $\mathbb{B}_\kappa$ with properties
$(1_\kappa)-(3_\kappa)$. When $|G|>\kappa$ then we can take any
subgroup $H\subset G$ of size $|H|=\kappa$ and find a subset
$\mathbb{B}_\kappa$ of $H$ with properties $(1_\kappa)-(3_\kappa)$
in $H$. Then the subset $\mathbb{B}_\kappa$ will have the
properties $(1_\kappa)-(3_\kappa)$ in the whole group.

So it remains to prove that such a set $\mathbb{B}_\kappa$ exists
in each
 group $G$ of size $\kappa$.

First we describe a sequence of  symmetric  subsets
$F_\alpha\subset G$ of size $\alpha$ such that
  $G=\bigcup_{\alpha < \kappa }F_{\alpha}$ and $F_{\alpha} \supset
  \bigcup_{\beta<\alpha}F_{\beta}$. Enumerate
  the group $G$ so that $G=\{g_\alpha : \alpha< \kappa \}$ and $g_0=e$. Then put
  $F_\alpha =\{g_\beta,-g_\beta : \beta< \alpha\}$ for all $\alpha<\kappa$.

We put
$$\mathbb{B}_\kappa=\bigcup_{ \alpha<\kappa } B_\alpha- B_\alpha$$
where a set $B_\alpha = \{b_\alpha^\beta:\beta<\alpha \}\subset G $ of
size $\alpha$  will be chosen later.

 To simplify notation we write $
 \mathbb{B}_{< \alpha}$ instead of $\bigcup_{ \beta<\alpha}
(B_\beta-B_\beta)$ and $
 \mathbb{B}_{> \alpha}$ instead of $\bigcup_{
\alpha\le\beta<\kappa}(B_\beta-B_\beta)$.
  By $B_\alpha^{<\beta}$ we shall denote the initial interval $\{b_\alpha^\gamma  : \gamma<\beta \}$ of $B_\alpha$.

Now we are in a position to define a sequence of sets $B_\alpha $
forcing the set  $\mathbb{B}_\kappa$ to satisfy the properties
$(2_\kappa)$ and $(3_\kappa)$.
To ensure property $(3_\kappa)$ we  will also construct a transfinite sequence
of points $(h_\alpha)_{ \alpha<\kappa}$ of $G$ such that $h_\alpha \notin
F_\alpha + \mathbb{B}_\kappa.$

We start putting $B_0= \{e\}$ and taking any non-zero point $h_0\in G$. Assume that for some ordinal $\alpha<\kappa$ the sets $B_\beta$ and the
points $h_\beta$, $\beta<\alpha$, have been constructed.
Then pick any point $h_\alpha\in G$ with  $$h_\alpha \notin F_\alpha +
\mathbb{B}_{<\alpha}.$$ Such a point exists because the size of
the set $F_\alpha +\mathbb{B}_{<\alpha }$ is equals $\alpha<\kappa=|G|$. Let
$$H_\alpha=\{h_\beta,-h_\beta :  \beta\le \alpha \}.$$

Next we define inductively  elements of $B_\alpha=\{b_\alpha^\beta
: \beta<\alpha \}.$

We pick any $b_\alpha ^0$ with  $b_\alpha^0\in G\backslash
\mathbb{B}_{<\alpha}$. Next, we chose $b_\alpha^\beta$  with

\begin{enumerate}
\item [(a)] $b_\alpha^\beta \notin
B_\alpha^{<\beta}+F_\alpha+\mathbb{B}_{<\alpha}$ ;
\item [(b)] $b_\alpha^\beta \notin B_\alpha^{<\beta}
-B_\alpha^{<\beta}+ B_\alpha^{<\beta}+ F_\alpha $;
\item [(c)]$ b_\alpha^\beta  \notin B_\alpha^{<\beta}+ F_\alpha
+H_\alpha
 $.
\end{enumerate}

To ensure properties (a),(b),(c) we have to avoid the  sets of
size $\alpha$, which is possible because $|G|=\kappa$. 

Now let us
prove that the constructed set $\mathbb{B}_\kappa$ satisfies
the properties $(1_\kappa)-(3_\kappa)$. In fact, the property
$(1_\kappa)$ is evident while $(3_\kappa)$ follows immediately
from (c).
It remains to prove
\smallskip

\textbf{Claim.} {\em The set $\mathbb{B}_\kappa$ has property
$(2_\kappa).$}
\smallskip

Let $B_\kappa$ be a subset of $ G$ of size $|B_\kappa|=\kappa$. Fix any
pairwise distinct points $c_1,c_2,c_3\in B_\kappa$.

If $B_\kappa-B_\kappa\subset \mathbb{B}_\kappa$ then
$B_\kappa\subset \bigcap_{i=1}^3 (c_i+\mathbb{B}_\kappa)$ and 
$\kappa=|B_\kappa|\leq | \bigcap_{i=1}^3 (c_i+\mathbb{B}_\kappa)|$.

So to prove our claim it is enough to show that $ |\bigcap_{i=1}^3
(c_i+\mathbb{B}_\kappa )|<\kappa.$ Find an ordinal $\alpha<\kappa$
such that $c_p-c_q\in F_\alpha$ for any $1 \leq p,q \leq 3$.
Assuming that  $ |\bigcap_{i=1}^3 (c_i+\mathbb{B}_\kappa) |=\kappa$
we may find a point $b\in \bigcap_{i=1}^3
(c_i+\mathbb{B}_{>\alpha})\backslash \{c_i\} $. A contradiction
will be reached in three steps.
\smallskip

\textbf{Step 1.} {\em First  show that there is $ \beta
>\alpha$ with $b\in \bigcap_{i=1}^3(c_i+B_\beta).$}
\smallskip

 Otherwise, $b-c_p\in B_\gamma-B_\gamma$
and $b-c_q\in B_\beta - B_\beta$ for some  $\gamma > \beta
>\alpha $ and some $p\ne q $. Find $i,j<\gamma$ with $b-c_p
=b_\gamma^i-b_\gamma^j.$ The inequality $b\ne c_p$  implies $i\ne
j$.

If $i < j$ then $b_\gamma^j= b_\gamma^i-b+c_p =b_\gamma^i-b+c_q-c_q+c_p\subset
b_\gamma^i-B_\beta+B_\beta+F_\gamma
\subset  B_\gamma^{<j} +\mathbb{B}_{<\gamma}+ F_\gamma,$
which contradicts (a).

If $i > j$ then $b_\gamma^i= b_\gamma^j+b-c_p = b_\gamma^j+B_\beta-B_\beta
+c_q-c_p \subset  B_\gamma^{<i} +\mathbb{B}_{<\gamma}+ F_\gamma,$ which again
contradicts (a).
\smallskip

\textbf{Step 2.} {\em We claim that if $b-c_p=b_\beta^i-b_\beta^j$
and $b-c_q=b_\beta^s-b_\beta^t$ then $\max \{i,j\}= \max
\{s,t\}.$}
\smallskip

It follows from the hypothesis that $c_q-c_p
=b_\beta^i-b_\beta^j+b_\beta^t-b_\beta^s.$
To obtain a contradiction assume that $\max \{i,j\}> \max
\{s,t\}.$ 

If $j < i$ then 
$b_\beta^i=c_q-c_p+b_\beta^j-b_\beta^t+b_\beta^s\in F_\beta +
B_\beta^{<i} -B_\beta^{<i}+B_\beta^{<i},$
which contradicts (b).

If $i < j$ then 
$b_\beta^j=c_p-c_q+b_\beta^i+b_\beta^t-b_\beta^s\in F_\beta +
B_\beta^{<j} +B_\beta^{<j}-B_\beta^{<j},$
again  a contradiction with (b).

\smallskip
\textbf{Step 3.} According to the previous step  there exists $
\beta>\alpha $ and $l$ such that
\smallskip

$b-c_1= b_\beta^i-b_\beta^j$ where $\max\{i,j\}$ is equal to $l$;

$b-c_2= b_\beta^s-b_\beta^t$ where  $\max\{s,t\}$ is equal to $l$;

$b-c_3= b_\beta^q-b_\beta^r$ where  $\max\{q,r\}$ is equal to $l$.

In this case we obtain a dichotomy: either among three numbers
$i,s,q$ two are equal to $l$ or among $j,t,r$ two are equal to
$l$. 
In the fist case we lose no generality assuming that $i=s=l;$ in
the second, that $j=t=l$.

In the first case we get 
$ F_\alpha \ni c_2-c_1= b_\beta^t-b_\beta^j,$ which contradicts (a).
In the second case we get 
$ F_\alpha \ni c_2-c_1= b_\beta^i-b_\beta^s,$ 
which contradicts (a) again.

Therefore, there is no $b\in \bigcap_{i=1}^3
(c_i+\mathbb{B}_{>\alpha})\backslash \{c_i\} $ and hence $|
\bigcap_{i=1}^3 (c_i+\mathbb{B}_{>\alpha})|<\kappa. $
\end{proof}

 {\em Acknowledgements.} The author
expresses her sincere thanks to Taras Banakh and Igor Protasov  for valuable and
stimulating discussions on the subject of the paper.


\begin{thebibliography}{}
\bibitem[BL$_1$]{BL1} T.~Banakh, N.~Lyaskovska, {\em Weakly P-small not
P-small subsets in Abelian groups}, Algebra and Discrete Mathematics, N.3 (2006), 29-34. 

\bibitem[BL$_2$]{BL2} T.~Banakh, N.~Lyaskovska, {\em Weakly P-small not
P-small subsets in groups}, Intern. J. of Algebra and Computations, {\bf 18}:1 (2008), 1-6.  

\bibitem[CS]{CS} J.~Convey, N.~Sloane, {\em Sphere packings, lattices, and groups,}  Springer, 1993.

\bibitem[DP]{DP} D.~Dikranjan, I.~Protasov, {\em Every infinite group can be
generated by a P-small subset}, General Applied Topology, {\bf 7}:2 (2006), 265-268.

\bibitem[Fu]{Fu} L.Fuchs, {\em Infinite abelian groups}, Academic Press, NY, 1970.


\end{thebibliography}
\end{document}